\documentclass[11pt,a4paper,draft]{amsproc}
\usepackage{amssymb}
\usepackage[hyphens]{url} 
\usepackage[dvips]{graphicx}
\theoremstyle{plain}
 \newtheorem{thm}{Theorem}[section]
 
 \newtheorem{lem}{Lemma}[section]
 
\theoremstyle{definition}

\theoremstyle{remark}
\usepackage{xcolor}
 
 \numberwithin{equation}{section}
\renewcommand{\leq}{\leqslant}
\renewcommand{\geq}{\geqslant}

\usepackage{geometry}
\newcommand{\<}{\langle}
\renewcommand{\>}{\rangle}
\renewcommand{\d}{\textrm{d}}

\usepackage{lipsum}
\newcommand\blfootnote[1]{%
	\begingroup
	\renewcommand\thefootnote{}\footnote{#1}%
	\addtocounter{footnote}{-1}%
	\endgroup
}

\title{AN ITERATIVE METHOD FOR SOLVING ELLIPTIC BVP IN ONE-DIMENSION}

\author{\bfseries Christian O. Bernal Zelaya$^{1}$, ~Prosper Torsu$^{1}$\\} 

\begin{document}
{\begin{flushleft}\baselineskip9pt\scriptsize
\end{flushleft}}
\vspace{18mm} \setcounter{page}{1} \thispagestyle{empty}

%The method  provides an efficient and computationally less expensive framework for solving such class of problems.

\maketitle
\vspace{-.5in}
\blfootnote{$^1$Department of Mathematics, California State University, Bakersfield, CA, U.S.A.}
\begin{abstract}
This paper presents a decomposition method for solving elliptic boundary value problems in one-dimension. The method is  an improvement to an existing technique for approximating elliptic systems. It is demonstrated to be computationally superior to the original formulation as less computations are required to obtain an approximation of the same accuracy. Convergence of the method is justified and supported by some theoretical results. We show that for a sufficiently smooth forcing data, the method always converge for a relatively small truncation order. The method is tested using some problems with exact solutions.
\end{abstract}

\section{Introduction} 
We consider the one-dimensional elliptic boundary-value problem (BVP) below:
\begin{equation}\label{eq:bvp}
\begin{cases}
\begin{aligned}
-\frac{\d}{\d x}\left(\kappa(x) \frac{\d u}{\d x} \right) &= f(x),~~~x\in (0, L)\\ 
u(0) = \alpha,~~~&-\kappa(L)\frac{\d u}{\d x}(L) = \beta.
 \end{aligned}
\end{cases}
\end{equation}
where $\kappa:(0,L) \rightarrow \mathbb{R}^+$ with $0 < \kappa_{\text{min}} \leq \kappa(x) \leq \kappa_{\text{max}}$ for all $x \in (0,L)$. Also, $f \in L^2(0,L)$, $\alpha$ and $\beta$ are given constants. Equation \eqref{eq:bvp} appears in the modeling of many physical events in physics, mechanics, and engineering. For instance, this model appears in the modeling of underground fluid flow, heat transfer, as well as applications in astronomy and electromagnetism.  In more recent years, Gou and Zhou have been able to figure out the exact traveling wave solutions of nonlinear Partial equations.  One of the more important roles in the studying of nonlinear physical phenomena is solving wave solutions of nonlinear PDE's, and recently these exact solutions help find these new phenomena.  An example of this can be seen in the waves of the ocean. The long ocean waves that are created by an island, can be modeled by the elliptic boundary-value problem Helmholtz equation \cite{krutitskii2001oblique}

Elliptic equations in simple forms can be solved using separation of variables in which closed form solutions can be obtained. Even though this a desirable, the method of separable equations does not always produce good approximations. For instance, if the magnitude of $a$ is very small, \eqref{eq:bvp} becomes a stiff problem where a series approximation is highly inaccurate because the Taylor expansion around zero may not converge. More importantly, we can obtain analytical solutions only for certain forms of $a$ and $f$. Such is why numerical methods have played a significant role in understanding solutions to differential equations. The three commonly used numerical methods are the Finite Volume Method (FVM), Finite Element Method (FEM) and, Finite Difference (FD). The Finite Volume Method \cite{babuska2001finite,chen2012higher,ewing2002accuracy,eymard2000finite} uses the divergence theorem to convert partial differential equations containing a divergence term into integral equations. The resulting equations are solved at the interface of neighboring elements. Finite Element Method \cite{jin2015finite,reddy2019introduction,zieli1992introduction,zienkiewicz1977finite, zienkiewicz2005finite} usually used in stress analysis, especially for determining displacements and stress within a structure. The method uses calculus of variations to convert \eqref{eq:bvp} into a minimization problem. Firstly, \eqref{eq:bvp} is transformed into a weak form. To obtain the discrete problem, the region of integration is divided into smaller non-overlapping subdomains, called elements, on which the problem is approximated using basis functions. The Finite Difference method \cite{wong2011exact} is very simple and straightforward to implement. It approximates differential equations by replacing derivatives by a corresponding difference quotient. In particular, the value of the unknown at discrete points is estimated by solving a finite number of equations involving the nodal values of nearby points.

This study is inspired by \cite{lu2007explicit}. In this article, Lu and Zhang studied a decomposition method where the unknown is expressed as an infinite series, which results in a collection of simpler problems to be solved. From a computational point of view, the method is limited by an inevitable repeated inversion of the Laplacian. In particular, the accuracy of the approximation gets better by increasing the number of terms in the partial sum. But, this leads to more inversions of the Laplacian. In this study, we derive an equivalent method which is more efficient and an approximation of the same accuracy. In contrast to the method proposed by Lu and Zhang where the Laplacian is inverted multiple times, the new formulation presented in this study inverts the Laplacian just two times. This means that we can increase the accuracy of our approximation without increasing the number of matrix inversions.

The rest of the paper is organized as follows: In Section \ref{sec:Definitions}, we present the basic definitions and notations that are used throughout the paper. Section \ref{sec:modelp} presents the mathematical model and its variational formulation. In Section \ref{sec:originalmethod} is an overview of the method studied by Lu and Zhang applied to \ref{sec:modelp}. In Section \ref{sec:impovedmethod}, we proposed an improved method and support it with numerical examples in Section  \ref{sec:examples}. We also provide a formal analysis of the method in Section \ref{sec:analysis} by showing that the approximation obtained using the proposed method converges to the true solution. Section \ref{sec:comclusion} concludes this study with an overview of proposed method and how it compares to the original method.

\section{Definitions and Notations.}\label{sec:Definitions}

Let $(0, L)$ be an interval of the real line. We define $L^2(0, L)$ as the space of square-integrable functions over $(0, L)$. It is endowed with the scalar product 
\begin{equation} \label{eq:ip}
\<u, v\> =  \int_{0}^L u(x)v(x)\d x. 
\end{equation}
$H^m(0, L)$ is a Hilbert space which contain functions whose weak derivatives of order up to $m$ are square-integrable over $(0, L)$. Particularly, $m=1$ corresponds to the space of functions defined on $(0, L)$ which, together, with their first derivatives belong to $L^2(0, L)$. That is
\begin{equation} \label{eq:Hminnn}
H^1(0, L) = \left\{v \in L^2(0,L): \textstyle \int_0^L (|v'|^2 + v^2) \d x < \infty \right\}.
\end{equation}
This space is equipped with scalar product defined by 
\begin{equation} \label{eq:Hminn}
\<u, v\>_{H^1} = \int_0^L \left(u'v' + uv\right)\, \d x.
\end{equation}
which induces the norm
\begin{equation} \label{eq:Hminn}
\|v\|_{H^1} =  \sqrt{\<v,v\>_{H^1}},
\end{equation}
where $\|\cdot \|$ is the usual $L^2$ - norm on $(0, L)$. For functions that vanish at $x = 0$ and $x = L$, the subspace $H^1_0(0, L)$ of $H^1(0, L)$ is used. In explicitly form,
\begin{equation} \label{eq:Hminnn}
H^1_0(0, L) = \left\{v \in H^1(0, L) : v(0) = v(L) = 0 \right\}.
\end{equation}
In the next section present the variational problem and some notations which are used in later sections.

\section{Model problem and variational formulation.}\label{sec:modelp}
 We start by defining the space $V = \left\{ v \in H^1(0,L): v(0) = 0 \right\}$. Then, the variational formulation for \eqref{eq:bvp} is to find $u \in V$ such that 
\begin{equation}\label{eq:fv}
a(u,v) = \ell(v) - \beta v(L),~\text{for all } v\in V.
\end{equation}
The bilinear form $a$ and the linear functional $\ell$ are defined respectively as
\begin{equation}\label{eq:bilinearform}
a(u,v) = \int_{0}^L \kappa(x)u'(x)v'(x)\,\d x, 
\end{equation}
and
\begin{equation}\label{eq:linearfunc}
\ell(v) = \int_{0}^L f(x)v(x)\,\d x.
\end{equation}
According to Lax-Milgram theorem, \eqref{eq:fv} has a unique solution in $V$.

\subsection{Existing work: An approximation}\label{sec:originalmethod}
The following approach is a direct consequence of the study in \cite{lu2007explicit}. We start by assuming that $u \in V$ can be written in the form 
\begin{equation}\label{eq:udecompu}
u(x) = \sum_{j=0}^\infty u_j(x) 
\end{equation}
Since $\kappa$ is assumed to be positive, we can find a smooth function $\psi$ such that $\psi = \log \kappa$. This means that 
\begin{equation}\label{eq:phi}
\kappa = \exp(\psi) = \sum_{j=0}^\infty \frac{\psi^j }{j!}
\end{equation}
Substituting \eqref{eq:udecompu} and \eqref{eq:phi} into the differential equation in \eqref{eq:bvp} gives
$$
-\frac{\d}{\d x}\left[\left(1 + \psi + \frac{1}{2!}\psi^2 + \frac{1}{3!}\psi^3 + \cdots \right) \left(\frac{\d u_0}{\d x} + \frac{\d u_1}{\d x} + \frac{\d u_2}{\d x} + \frac{\d u_3}{\d x} + \cdots\right)\right] = f(x)
$$
Upon expansion and matching the terms so that the subscript of $u$ and exponent of $\psi$ sum up to the same integer, we have the following equations:
\begin{equation}\label{eq:serieseq1}
\begin{split}
-\frac{\d^2u_0}{\d x^2} &= f(x),\\[1ex]
-\frac{\d^2u_1}{\d x^2} &- \frac{\d}{\d x}\left(\psi \frac{\d u_0}{\d x} \right) = 0,\\[1ex]
-\frac{\d^2u_2}{\d x^2} &- \frac{\d}{\d x}\left(\psi \frac{\d u_1}{\d x} \right) + \frac{1}{2!}\frac{\d}{\d x}\left(\psi^2 \frac{\d u_0}{\d x} \right) = 0,\\[1ex]
\vdots\\
-\frac{\d^2u_n}{\d x^2} &- \sum_{j = 1}^n \frac{1}{j!} \frac{\d}{\d x}\left(\psi^j \frac{\d u_{n-j}}{\d x} \right) = 0.
\end{split}
\end{equation}
and so on. We derive the boundary conditions similarly. That is 
\begin{equation}\label{eq:bcdserieseq}
u_0(0) = \alpha,~~~u_1(0) = 0,~~~u_2(0) = 0, ~~\cdots~~~u_M(0) = 0,\cdots
\end{equation}
and
\begin{equation}\label{eq:bcnserieseq}
\begin{split}
-\frac{\d u_0}{\d x}(L) &= \beta\\[1ex]
-\frac{\d u_1}{\d x}(L) &- \psi(L) \frac{\d u_0}{\d x}(L) = 0,\\[1ex]
-\frac{\d u_2}{\d x}(L) &- \psi(L) \frac{\d u_1}{\d x}(L) + \frac{\psi^2(L)}{2!} \frac{\d u_0}{\d x}(L) = 0,\\
\vdots\\
-\frac{\d u_n}{\d x}(L) &- \sum_{j = 1}^n \frac{\psi^j(L)}{j!} \frac{\d u_{n-j}}{\d x}(L) = 0
\end{split}
\end{equation}
and so on. By collecting the differential equations in \eqref{eq:serieseq1} and the corresponding boundary conditions in \eqref{eq:bcdserieseq} and \eqref{eq:bcnserieseq}, we get the following boundary value problems.
\begin{equation}\label{eq:bvps}
\begin{aligned}
&
\begin{cases}
\begin{aligned}
-\frac{\d^2u_0}{\d x^2} &= f(x),\\[1ex]
u_0(0) &= \alpha,~~~-\frac{\d u_0}{\d x}(L) = \beta.
 \end{aligned}
\end{cases}\\[1ex]
&
\begin{cases}
\begin{aligned}
-\frac{\d^2u_1}{\d x^2} & = \frac{\d}{\d x}\left(\psi \frac{\d u_0}{\d x} \right),\\[1ex]
u_1(0) &= 0,~~~-\frac{\d u_1}{\d x}(L) - \psi(L) \frac{\d u_0}{\d x}(L) = 0.
 \end{aligned}
\end{cases}\\[1ex]
&
\begin{cases}
\begin{aligned}
-\frac{\d^2u_2}{\d x^2} & = \frac{\d}{\d x}\left(\psi \frac{\d u_1}{\d x} \right) + \frac{1}{2!}\frac{\d}{\d x}\left(\psi^2 \frac{\d u_0}{\d x} \right),\\[1ex]
u_2(0) &= 0,~~~- \psi(L) \frac{\d u_1}{\d x}(L) - \frac{\psi^2(L)}{2!} \frac{\d u_0}{\d x}(L) = 0.
 \end{aligned}
\end{cases}\\[1ex]
&
\hspace*{1.5in}\vdots\\
&
\begin{cases}
\begin{aligned}
-\frac{\d^2u_m}{\d x^2} & = \sum_{j = 1}^m \frac{1}{j!} \frac{\d}{\d x}\left(\psi^j \frac{\d u_{n-j}}{\d x} \right),\\[1ex]
u_m(0) &= 0,~~~-\frac{\d u_m}{\d x}(L) - \sum_{j = 1}^m \frac{\psi^j(L)}{j!} \frac{\d u_{m-j}}{\d x}(L) = 0
 \end{aligned}
\end{cases}
 \end{aligned}
\end{equation}
The first $m$ equations are then solved for and, an approximation to \eqref{eq:bvp} is obtained using
\begin{equation}\label{eq:uapprox}
U_m = \sum_{j=0}^m u_j
\end{equation}
For the purpose of analysis and discussions to follow, we shall denote
\begin{equation}\label{eq:newvf}
a_j(u,v) = \frac{1}{j!}\int_0^L\psi^j(x)u'(x)v'(x)\, \d x.
\end{equation}
Then, the variational formulation for the method in \eqref{eq:bvps} is 
\begin{equation}\label{eq:oldmethod}
\begin{cases}
\begin{split}
&\text{Find $u_0 \in V$ satisfying } a_0(u_0, v) = \ell(v) - \beta v(L),~  ~~\text{forall } v \text{ in }  V\\
&\text{Find $u_m \in V$ satisfying } {a}_0({u}_m, v) +\sum_{j=1}^{m} a_j({u}_{m-j}, {v}) = 0,~  ~~\text{forall } v \text{ in }  V
\end{split}
\end{cases}
\end{equation}
The drawback of this method is that although it allows us to rewrite the original problem as a collection of simpler and identical problem, the process of obtaining the solutions is sequential. This means that, if we are interested in an approximation of order $M$, we must invert the same Poisson operator $M$ times in order to provide the source data for $u_M$. As a result, we must perform $M+1$ matrix inversions to obtain $u_M$, and this can make the method very expensive to solve. In the following section, we present an improve method that reduces the amount of matrix inversion from $M+1$ to just $2$, regardless of $M$. More importantly, the accuracy of the approximation obtained using $M+1$ inversions is maintained.

\subsection{Undoing the recursive relations}\label{sec:impovedmethod}
 As pointed out above, in order to solve for $u_M$, we must first solve for $u_0, u_1, \cdots, u_{M-1}$ which is more expensive than solving the original problem. We are going to undo the dependence of higher order terms on all lower order once and show that one can obtain $u_M$ directly by solving for $u_0$ and $u_M$. This way, we avoid calculations of $u_1, \cdots, u_{M-1}$ which is inevitable in the original method. We assume that $u_0$ is available, which means that solving the BVP for $u_0$ is unavoidable. We show the derivation for the differential equations only because the boundary conditions follow trivially.

\subsubsection{Equation for $u_1$}
The equation for $u_1$ is 
\begin{equation*}
-\frac{\d^2u_1}{\d x^2} = \frac{\d}{\d x}\left(\psi \frac{\d u_0}{\d x} \right),
\end{equation*}
which involves $u_0$ only. This means that we can obtain $U_1 = u_0 + u_1$ by solve two BVPs- for $u_0$ and $u_1$. Now, by integrating the equation above, we have 
$$
-\int_x^L\frac{\d^2u_1}{\d s^2}\, ds - \int_x^L\frac{\d}{\d s}\left(\psi \frac{\d u_0}{\d s} \right)\, ds = 0
$$
It follows from the Fundamental Theorem of Calculus (FTC) that
$$
-\frac{\d u_1}{\d x}(L) + \frac{\d u_1}{\d x} - \psi(L) \frac{\d u_0}{\d x}(L) + \psi \frac{\d u_0}{\d x} = 0
$$
Using the boundary condition for $u_1$, namely $-\frac{\d u_1}{\d x}(L) - \psi(L) \frac{\d u_0}{\d x}(L) = 0$, we have 
\begin{equation}\label{eq:du1}
\frac{\d u_1}{\d x} + \psi \frac{\d u_0}{\d x} = 0~~\text{or}~~\frac{\d u_1}{\d x} = -\psi \frac{\d u_0}{\d x}
\end{equation}

\subsubsection{Equation for $u_2$}
Here, we write the differential equation for $u_2$ in terms of $u_0$ only by eliminating $u_1$. We recall that $u_2$ satisfies
\begin{equation}\label{eq:deu2}
-\frac{\d^2u_2}{\d x^2} = \frac{\d}{\d x}\left(\psi \frac{\d u_1}{\d x} \right) + \frac{1}{2!}\frac{\d}{\d x}\left(\psi^2 \frac{\d u_0}{\d x} \right) = 0
\end{equation}
Using \eqref{eq:du1}, this equation becomes
\begin{equation*}
\begin{split}
-\frac{\d^2u_2}{\d x^2} &= \frac{\d}{\d x}\left(-\psi^2 \frac{\d u_0}{\d x} \right) + \frac{1}{2!}\frac{\d}{\d x}\left(\psi^2 \frac{\d u_0}{\d x} \right)\\[1ex]
& = -\frac{\d}{\d x}\left(\frac{\psi^2}{2!} \frac{\d u_0}{\d x} \right)
\end{split}
\end{equation*}
Integration \eqref{eq:deu2} and using the FTC and \eqref{eq:du1}, we have 
\begin{equation}\label{eq:du2}
\frac{\d u_2}{\d x} = \frac{\psi^2}{2!} \frac{\d u_0}{\d x}
\end{equation}

\subsubsection{Equation for $u_3$}
We start with the differential equation for $u_3$, which is 
$$
-\frac{\d^2u_3}{\d x^2} = \frac{\d}{\d x}\left(\psi \frac{\d u_2}{\d x} \right) + \frac{1}{2!}\frac{\d}{\d x}\left(\psi^2 \frac{\d u_1}{\d x} \right) + \frac{1}{3!}\frac{\d}{\d x}\left(\psi^3 \frac{\d u_0}{\d x} \right)
$$
Using \eqref{eq:du1} and \eqref{eq:du2}
\begin{equation*}
\begin{split}
-\frac{\d^2u_3}{\d x^2} &= \frac{1}{2!}\frac{\d}{\d x}\left(\psi \frac{\d u_2}{\d x} \right) + \frac{1}{2!}\frac{\d}{\d x}\left(\psi^2 \frac{\d u_1}{\d x} \right) + \frac{1}{3!}\frac{\d}{\d x}\left(\psi^3 \frac{\d u_0}{\d x} \right)\\[1ex]
& =  \frac{1}{2!}\frac{\d}{\d x}\left(\psi^3 \frac{\d u_0}{\d x} \right) + \frac{1}{2!}\frac{\d}{\d x}\left(-\psi^3 \frac{\d u_1}{\d x} \right) + \frac{1}{3!}\frac{\d}{\d x}\left(\psi^3 \frac{\d u_0}{\d x} \right)\\[1ex]
& = \frac{\d}{\d x}\left(\frac{\psi^3}{3!} \frac{\d u_0}{\d x} \right)
\end{split}
\end{equation*}
Successfully, $u_1$ and $u_2$ have been eliminated from the equation for $u_3$. By repeated application of the  FTC and the boundary conditions, we have 
\begin{equation}\label{eq:du2}
-\frac{\d^2u_j}{\d x^2} = 
\frac{(-1)^{j+1}}{j!}\frac{\d}{\d x}\left(\psi^j \frac{\d u_0}{\d x} \right)
\end{equation}
for $j = 1,2,3,\cdots,m$. This way, $u_1 + u_2 + \cdots + u_m$ satisfies
\begin{equation*}
-\frac{\d^2}{\d x^2}(u_1 + u_2 + \cdots + u_m) = \frac{\d}{\d x}\left[\left(\psi - \frac{\psi^2}{2!} + \frac{\psi^3}{3!} - \cdots + (-1)^{m+1}\frac{\psi^m}{m!} \right)\frac{\d u_0}{\d x} \right]
\end{equation*}
We observe that the coefficient on the right hand side is the partial sum of $e^{-\psi} - 1$. Therefore, by appending $u_0$, the approximation $U_m = u_0 + u_1 + u_2 + \cdots + u_m$ satisfies
\begin{equation}
-\frac{\d^2 U_m}{\d x^2} - \frac{\d}{\d x}\left[\left(1 - \psi + \frac{\psi^2}{2!} - \frac{\psi^3}{3!} + \cdots + (-1)^{m+1}\frac{\psi^m}{m!} \right)\frac{\d u_0}{\d x} \right] = 0
\end{equation}
or simply 
\begin{equation*}
-\frac{\d^2 U_m}{\d x^2} - \frac{\d}{\d x}\left(G_m(x)\frac{\d u_0}{\d x} \right) = 0
\end{equation*}
where 
\begin{equation}\label{eq:gm}
G_m = \displaystyle \sum_{j = 0}^m \frac{(-\psi)^{j}}{j!}
\end{equation}
Here, an approximation for the exact solution is $U_m$. Clearly, this approach eliminates the dependence of $u_m$ on all of $u_1, u_2,\cdots,u_{m-1}$. The difference between the new derivation and the original approach is that, for the old method we first solve for  $u_0, u_1, u_2,\cdots,u_{m-1}$ in order to create the create the source terms for $u_m$. With the new method, we simply need $u_0$ to create the source terms for $u_m$. We present the new method as the following variational formulation:
\begin{equation}\label{eq:newmethod}
\begin{cases}
\begin{split}
&\text{Find $u_0 \in V$ satisfying } a_0(u_0, v) = \ell(v) - \beta v(L),~~~\text{forall } v \text{ in }  V\\[1ex]
&\text{Find $u_M \in V$ satisfying } {a}_0({u}_M, v) + (G_Mu'_0, v') = 0,~~~\text{forall } v \text{ in }  V
\end{split}
\end{cases}
\end{equation}
where $G_m$ is defined in \eqref{eq:gm}.

\section{Numerical examples}\label{sec:examples}
Here, we present some numerical examples to validate the proposed method. For simplicity, all approximations are calculated on the domain $(0,1)$ with $\alpha = \beta = 0$. We remark that the total error associated with this method consists of two parts: the error associated with the decomposition and secondly, the error due to solving a discrete problem.

\subsection{Example 1} We set $\kappa(x) = 1+x^2,~f(x) = 1$, so that the exact solution is 
\begin{equation}
u(x) = \tan^{-1}(x) - \frac{1}{2}\ln(1 + x^2) 
\end{equation}
The following table shows the associated $L^2$ errors with various mesh widths:\\[-1ex]
\begin{center}
\begin{tabular}{ |c|c|c|c|c|c| } 
 \hline
 & $M=2$ & $M=4$ & $M=6$ & $M=8$ & $M=10$ \\\hline 
$N = 2^3 $ & $1.2839(-02)$ & $1.3365(-02)$ & $1.3384(-02)$ & $1.3384(-02)$ & $1.3384(-02)$ \\\hline 
$N = 2^5$ & $3.3582(-03)$ & $3.2583(-03)$ & $3.2636(-03)$ & $3.2637(-03)$ & $3.2637(-03)$ \\\hline 
$N = 2^7 $ & $1.1925(-03)$ & $8.1399(-04)$ & $8.1288(-04)$ & $8.189(-04)$ & $8.1289(-04)$ \\\hline 
$N = 2^9 $ & $6.9101(-04)$ & $2.0603(-04)$ & $2.0307(-04)$ & $2.0306(-04)$ & $2.0306(-04)$ \\\hline 
$N = 2^{11} $ & $5.7681(-04)$ & $5.4292(-04)$ & $5.0773(-05)$ & $5.0756(-05)$ & $5.0756(-05)$ \\\hline 
\end{tabular}
\end{center}

\subsection{Example 2} Here, we set $\kappa(x) = (1-0.5\sin(10\pi x))^{-1},~f(x) = 1$. The exact solution is  
\begin{equation}
u(x) = \frac{1}{200\pi^2} \left[\sin(10 \pi x) + 10 \pi (1- x) \cos(10 \pi x) + 100\pi^2 x(2- x)\right] - \frac{1}{200\pi}
\end{equation}

The $L^2$ errors associated with various mesh widths are shown below:\\[-1ex]
\begin{center}
\begin{tabular}{ |c|c|c|c|c|c| } 
 \hline
 & $M=2$ & $M=4$ & $M=6$ & $M=8$ & $M=10$ \\\hline 
$N = 2^3 $ & $4.7783(-02)$ & $4.5079(-02)$ & $4.4994(-02)$ & $4.993(-02)$ & $4.993(-02)$ \\\hline 
$N = 2^5$ & $2.1835(-02)$ & $1.9095(-02)$ & $1.9029(-02)$ & $1.9029(-02)$ & $1.9029(-02)$ \\\hline 
$N = 2^7 $ & $5.6535(-03)$ & $2.9191(-03)$ & $2.8536(-03)$ & $2.8529(-03)$ & $2.8529(-03)$ \\\hline 
$N = 2^9 $ & $3.4024(-03)$ & $6.5572(-04)$ & $5.9061(-04)$ & $5.8993(-04)$ & $5.8992(-04)$ \\\hline 
$N = 2^{11} $ & $2.9603(-03)$ & $2.0610(-04)$ & $1.4076(-04)$ & $1.4007(-04)$ & $1.4007(-04)$ \\\hline 
\end{tabular}
\end{center}

\subsection{Example 3} Here, we set $\kappa(x) = (x + 1)^2,~f(x) = x/(x+1)$. The exact solution is  
\begin{equation}
u(x) = \frac{(3 - \ln 2)x - (2 + x) \ln(1+x)}{1 + x}
\end{equation}

The $L^2$ errors associated with various mesh widths are shown below: \\[-1ex]
\begin{center}
\begin{tabular}{ |c|c|c|c|c|c| } 
 \hline
 & $M=2$ & $M=4$ & $M=6$ & $M=8$ & $M=10$ \\\hline 
$N = 2^3 $ & $3.5860(-03)$ & $4.1612(-03)$ & $4.5709(-03)$ & $4.5944(-03)$ & $4.5950(-03)$ \\\hline 
$N = 2^5$ & $4.6503(-03)$ & $1.0192(-03)$ & $1.0143(-03)$ & $1.0190(-03)$ & $1.0191(-03)$ \\\hline 
$N = 2^7 $ & $5.2928(-03)$ & $4.1726(-04)$ & $2.4933(-04)$ & $2.4732(-04)$ & $2.4732(-04)$ \\\hline 
$N = 2^9 $ & $5.4764(-03)$ & $2.9496(-04)$ & $6.5412(-05)$ & $6.1425(-05)$ & $6.1373(-05)$ \\\hline 
$N = 2^{11} $ & $5.5237(-03)$ & $2.6897(-04)$ & $2.0318(-05)$ & $1.5381(-05)$ & $1.5315(-05)$ \\\hline 
\end{tabular}
\end{center}

\subsection{Example 4} We set $\kappa(x) = x^4 + e^{-x},~f(x) = -2\cos(\pi x)$. The semi-analytical solution is 
\begin{equation}
u(x) = \int_0^1\frac{2\sin \pi s}{\pi(s^4 + e^{-s})}ds - \int_x^1\frac{2\sin \pi s}{\pi(s^4 + e^{-s})}ds
\end{equation}
The $L^2$ errors associated with various mesh widths are calculated using $2^{15}$ elements for the reference solution. The results are presented in the table below:\\[-1ex]
\begin{center}
\begin{tabular}{ |c|c|c|c|c|c| } 
 \hline
 & $M=2$ & $M=4$ & $M=6$ & $M=8$ & $M=10$ \\\hline 
$N = 2^3$     & $9.1918(-02)$ & $9.3818(-02)$ & $9.3831(-02)$ & $9.3831(-02)$ & $9.3831(-02)$ \\\hline 
$N = 2^5$     & $1.5883(-02)$ & $1.7578(-02)$ & $1.7589(-02)$ & $1.7589(-02)$ & $1.7589(-02)$ \\\hline 
$N = 2^7$     & $2.4796(-03)$ & $4.0802(-03)$ & $4.0908(-03)$ & $4.0908(-03)$ & $4.0908(-03)$ \\\hline 
$N = 2^9$     & $7.0333(-04)$ & $9.7832(-04)$ & $9.8874(-04)$ & $9.8874(-04)$ & $9.8874(-04)$\\\hline 
$N = 2^{11}$ & $1.0989(-04)$ & $1.2989(-04)$ & $1.2989(-04)$ & $1.2989(-04)$ & $1.2989(-04)$\\\hline 
\end{tabular}
\end{center}

\section{An analysis}\label{sec:analysis}
In this section, we shall prove that the individual terms in the expansion of the solution are bounded. We will also prove that as the number of terms in the expansion increases, the approximations get more accurate, and finally converges to the exact solution. We start with the following preliminary results will be used.

\begin{lem}\label{eq:lemma1}
Given $v \in V$,  we can find $w \in V$ such that $w'(x) = \psi(x) v'(x)$, where $\psi \in L^2(0,L)$.

\begin{proof}
By taking $w(x) = \int_{0}^x \psi(s)v'(s)ds$, follows that $w'(x) = \psi(x)v'(x)$. Moreover, $w(0) = 0$.
\end{proof}
\end{lem}

\begin{lem} \label{lem:longsum}
Let $\psi \in L^2(0, 1)$. Then for any positive integer $M$, it holds that
$$
\sum_{j=M+1}^{\infty} \frac{\|\psi \|^j_{\infty}}{j!}  <  \frac{\|\psi\|_{\infty}^{M+1}}{(M+1)!} \exp (\|\psi\|_{\infty}).
$$
\begin{proof}
Direct manipulation shows that 
\begin{align*}
  \sum_{j=M+1}^{\infty} \frac{\|\psi\|_{\infty}^j}{j!} & = \frac{\|\psi\|_{\infty}^{M+1}}{(M+1)!} \left(1 + \frac{\|\psi\|_{\infty}}{(M+2)} + \frac{\|\psi\|_{\infty}^2}{(M+2)(M+3)} + \cdots \right)\\ 
  & <  \frac{\|\psi\|_{\infty}^{M+1}}{(M+1)!} \left(1 + \|\psi\|_{\infty} + \frac{\|\psi\|_{\infty}^2}{2!} + \frac{\|\psi\|_{\infty}^3}{3!} + \cdots \right) \\
 & =  \frac{\|\psi\|_{\infty}^{M+1}}{(M+1)!}\,\text{exp}(\|\psi\|_{\infty})
 \end{align*} 
 \end{proof}
\end{lem}

\begin{lem}\label{eq:lemma2}
Suppose that $\sum_{j=0}^{m} a_j({u}_{m-j}, {v}) = 0$. Then $\sum_{j = 1}^{m} ja_j(u_{m-j}, v) = 0$ for all $n< m$.
\begin{proof}
Assume that $\sum_{j = 0}^{n} a_j(u_{n-j}, v) = 0,~ ~~\text{for all}~~ v\in V$. Then
\begin{equation*}
\begin{split}
\sum_{j = 1}^m ja_{j}(u_{m-j}, v) &= a_{1}(u_{m-1}, v) + 2a_{2}(u_{m-2}, v) + \cdots+ ma_{m}(u_{0}, v)\\
 &= \<\psi u'_{m-1}, v'\> + 2 \cdot \frac{1}{2!} \<\psi^2 u'_{m-2}, v'\> + \cdots + m \cdot \frac{1}{m!} \<\psi^m u'_{0}, v'\> \\
 &= \<u'_{m-1}, \psi v'\> + \<\psi u'_{m-2}, \psi v'\> + \cdots + \frac{1}{(m-1)!} \<\psi^{m-1} u'_{0}, \psi v'\> \\
\end{split}
\end{equation*}
By Lemma \ref{eq:lemma4}, there exists  $w \in V$ such that $w' = \psi v'$. Therefore,

\begin{equation*}
\begin{split}
\sum_{j = 1}^m ja_{j}(u_{m-j}, v) &= \<u'_{m-1}, \psi v'\> + \<\psi u'_{m-2}, \psi v'\> + \cdots + \frac{1}{(m-1)!} \<\psi^{m-1} u'_{0}, \psi v'\> \\
&= \<u'_{m-1}, w'\> + \<\psi u'_{m-2}, w'\> + \cdots + \frac{1}{(m-1)!} \<\psi^{m-1} u'_{0}, w'\> \\
& = \sum_{j = 0}^{m-1} a_j(u_{m-1-j}, w)\\
& = 0 
\end{split}
\end{equation*}
A more general result is 
\begin{equation}
\<u'_{n} +  \psi u'_{n-1} + \tfrac{1}{2!} \psi^2 u'_{n-2}  + \cdots + \tfrac{1}{n!}\psi^n u'_{0}, ~\psi^q w' \> = 0,~~ ~~\text{for all}~~ q\geq 2.
\end{equation}
\end{proof}
\end{lem}

\begin{lem}\label{eq:lemma3}
Let $u_m$ satisfy the hypothesis of Lemma \ref{eq:lemma1}. Then 
$$
a_0(u_m, v) = (-1)^ma_m(u_0, v),~~~\text{for all},~~ v \in V~~ \text{ and }~~  m\geq 1
$$
\begin{proof}
This is a direct consequence of \eqref{eq:du2}.
%
%We induct on $m$. Clearly, the result holds for $m=1$; namely, 
%$$
%a_0(u_1, v) = - a_1(u_0, v) = (-1)^1a_1(u_0, v) \ \forall v\in V
%$$
%For $m = 2$, we consider the equation for $u_2$ given by 
%$$
%a_0(u_2, v) + a_1(u_1, v) + a_2(u_0, v) = 0
%$$
%By adding and subtracting $2a_2(u_0, v)$, this can be written as 
%\begin{equation*}
%\begin{split}
%0 &= a_0(u_2, v) + a_1(u_1, v) +[ 2a_2(u_0, v) - 2a_2(u_0, v)] + a_2(u_0, v) \\
%& = a_0(u_2, v) + a_1(u_1, v) + 2a_2(u_0, v) - a_2(u_0, v)
%\end{split}
%\end{equation*}
%By Lemma \ref{eq:lemma1}, $ a_1(u_1, v) + 2a_2(u_0, v) = 0$. Therefore, the last equation becomes
%$$
%a_0(u_2, v) =  a_2(u_0, v) = (-1)^2a_2(u_0, v)
%$$
%For $m\geq 3$, the result 
%$$
%a_0(u_m, v) = (-1)^ma_m(u_0, v)
%$$
%follows similarly by adding $\sum_{j = 1}^m ja_j(u_{m-j}, v) - \sum_{j = 1}^m ja_j(u_{m-j}, v)$ to the variational formulation for $u_m$, and the using Lemma \ref{eq:lemma1}.
\end{proof}
\end{lem}
\begin{lem}\label{eq:lemma4}
Suppose that $\sum_{j=0}^{m} a_j({u}_{m-j}, {v}) = 0$. Then $|u_m|_{H^1} \leq \tilde{C}_m \|\psi \|^m_{\infty} \|f \|_\infty$ for all $m \geq 0$.
\begin{proof}
This result follows from Lemma \ref{eq:lemma3}. First, the variational equation for $u_0$ is $a_0(u_0,v) = \<f,v\>$. It follows from Poincare's inequality that $|u_0|_{H^1} \leq C_0\|f\|_\infty$
for some $C_0$ that does not depend on $u_0$. The identity
$$
|u_m|_{H^1} \leq C_m \|\psi\|^m_{\infty}|u_0|_{H^1} \leq \tilde{C}_m \|\psi\|^m_{\infty}\|f\|_\infty
$$
is establish via induction and using Lemma \ref{eq:lemma2}.
\end{proof}
\end{lem}
Lemma \ref{eq:lemma3} shows that every term in the expansion of the solution to the BVP is bounded. Therefore, increasing the number of terms in the approximation for solution guarantees a stable solution.

\begin{lem} \label{lem:vf}
Let $u$ be the solution \eqref{eq:bvp}, and $u_0, u_M$ satisfy \eqref{eq:newmethod}. Then 
$$
\<u' - u'_M, v'\> =   \sum_{j = M+1}^{\infty} \frac{(-1)^{j}}{j!}\<\psi^j u'_0, v'\>~~  \text{ for all } v\in V
$$
\begin{proof}
It follows from the variations equations that 
$$
\<\kappa u',  v'\> = \<f, v\> - \beta v(L) = \<u'_0, v'\>
$$ 
which is equivalent to $\<\kappa u',  v'\> =  \<u'_0, v'\>$. Since $\kappa$ is invertible, in particular, one can write
$$
\<\kappa u' - u'_0, v'\> = 0~~~~\text{ or }~~~~\<u' - \kappa^{-1}u'_0, \kappa v'\> = 0
$$ 
Lemma \ref{eq:lemma1} establishes the existence of $w \in V$ for which $w' = \kappa v'$. Therefore, 
$$
\<u' - \kappa^{-1}u'_0, w'\> = 0 ~~~~\text{ or }~~~~ \<u',  w'\> =  \<\kappa^{-1}u'_0, w'\>
$$
Thus,
\begin{align*}
\<u' - u'_M, v'\> &=  \<u', v'\> - \<u'_M, v'\>\\
&= \<\kappa^{-1}u'_0, v'\> - \<B_Mu'_0, v'\>\\
&= \<(\kappa^{-1} - B_M)u'_0, v' \big\>
\end{align*}
But then
\begin{align*}
\kappa^{-1} - B_M = \exp(-\psi) - \sum_{j = M+1}^{\infty} \frac{(-1)^{j}}{j!}\psi^j = \sum_{j = M+1}^{\infty} \frac{(-1)^{j}}{j!}\<\psi^j u'_0, v' \>
\end{align*}
It follows that 
\begin{align*}
\<u' - u'_M, v'\> =   \sum_{j = M+1}^{\infty} \frac{(-1)^{j}}{j!}\<\psi^j u'_0, v' \>
\end{align*}
\end{proof}
\end{lem}
The following major theorem establishes the convergence of the iterative method.
 \begin{thm}\label{theo:main}
 Let  $u_M$ be the approximate in \eqref{eq:newmethod}. Then $\displaystyle \lim_{M \to \infty} \big|u - u_M \big|_{H^1} = 0$.
 \begin{proof}
By Lemma \ref{lem:vf}, we have 
 \begin{align*}
\<u' - u'_M, v'\> &=   \sum_{j = M+1}^{\infty} \frac{(-1)^{j}}{j!} \<\psi^j u'_0, v' \>
\end{align*}
In particular
\begin{align*}
\<u' - u'_M, u' - u'_M\> &=   \sum_{j = M+1}^{\infty}\frac{(-1)^{j}}{j!} \<\psi^j u'_0, u' - u'_M \>
\end{align*}
By Cauchy-Schwarz and triangle inequality, we get 
\begin{align*}
\big|u - u_M \big|_{H^1} \leq \left(\sum_{j = M+1}^{\infty} \frac{\|\psi\|_{\infty}^j}{j!}\right) \big| u_0\big|_{H^1} 
\end{align*}
It follows that 
\begin{align*}
\big|u - u_M \big|_{H^1} &< \frac{\|\psi\|_{\infty}^{M+1}}{(M+1)!} \,\exp(\|\psi\|_{\infty}) \big| u_0\big|_{H^1} 
\end{align*}
As $\psi$ is bounded, we conclude that 
\begin{align*}
\lim_{M \to \infty} \big|u - u_M \big|_{H^1} = 0
\end{align*}
Thus, $u_M$ converges to $u$ in $H^1 - $ norm.
\end{proof}
\end{thm}

\section{Conclusion}\label{sec:comclusion}
We studied a decomposition method for approximating elliptic boundary value problems in one-dimension. The method is shown to be efficient and produces accurate approximations. We have shown that within some acceptable tolerance, the method con- verges rapidly with respect to the truncation order $M$. For sufficiently smooth coefficients, only a few terms are needed to obtain accurate approximations. For highly variable coefficients however, a relatively large number of terms may be required for an accurate approximation. The method has been analyzed and shown to produce an approximation that converges in $H^1$-norm to the true solution for sufficiently large $M$. The analysis also provides a general framework for studying convergence properties of the original formulation. Some experimental examples were provided to support the effectiveness of the method.

\bibliographystyle{siamplain}
\bibliography{references}

\begin{thebibliography}{10}

\bibitem{babuska2001finite}
{\sc I.~Babuska, T.~Strouboulis, I.~Babu{\v{s}}ka, J.~R. Whiteman, et~al.},
  {\em The finite element method and its reliability}, Oxford university press,
  2001.

\bibitem{chen2012higher}
{\sc Z.~Chen, J.~Wu, and Y.~Xu}, {\em Higher-order finite volume methods for
  elliptic boundary value problems}, Advances in Computational Mathematics, 37
  (2012), pp.~191--253.

\bibitem{ewing2002accuracy}
{\sc R.~E. Ewing, T.~Lin, and Y.~Lin}, {\em On the accuracy of the finite
  volume element method based on piecewise linear polynomials}, SIAM Journal on
  Numerical Analysis, 39 (2002), pp.~1865--1888.

\bibitem{eymard2000finite}
{\sc R.~Eymard, T.~Gallou{\"e}t, and R.~Herbin}, {\em Finite volume methods},
  Handbook of numerical analysis, 7 (2000), pp.~713--1018.

\bibitem{jin2015finite}
{\sc J.-M. Jin}, {\em The finite element method in electromagnetics}, John
  Wiley \& Sons, 2015.

\bibitem{krutitskii2001oblique}
{\sc P.~Krutitskii}, {\em The oblique derivative problem for the helmholtz
  equation and scattering tidal waves}, Proceedings of the Royal Society of
  London. Series A: Mathematical, Physical and Engineering Sciences, 457
  (2001), pp.~1735--1755.

\bibitem{lu2007explicit}
{\sc Z.~Lu, D.~Zhang, and B.~A. Robinson}, {\em Explicit analytical solutions
  for one-dimensional steady state flow in layered, heterogeneous unsaturated
  soils under random boundary conditions}, Water resources research, 43 (2007).

\bibitem{reddy2019introduction}
{\sc J.~N. Reddy}, {\em Introduction to the finite element method}, McGraw-Hill
  Education, 2019.

\bibitem{wong2011exact}
{\sc Y.~S. Wong and G.~Li}, {\em Exact finite difference schemes for solving
  helmholtz equation at any wavenumber}, International Journal of Numerical
  Analysis and Modeling, Series B, 2 (2011), pp.~91--108.

\bibitem{zieli1992introduction}
{\sc T.~G. Zieli}, {\em Introduction to the finite element method},  (1992).

\bibitem{zienkiewicz1977finite}
{\sc O.~C. Zienkiewicz, R.~L. Taylor, P.~Nithiarasu, and J.~Zhu}, {\em The
  finite element method}, vol.~3, McGraw-hill London, 1977.

\bibitem{zienkiewicz2005finite}
{\sc O.~C. Zienkiewicz, R.~L. Taylor, and J.~Z. Zhu}, {\em The finite element
  method: its basis and fundamentals}, Elsevier, 2005.

\end{thebibliography}
\end{document}